\documentclass[reqno]{amsart} 

\usepackage{amsfonts}
\usepackage{amsmath}
\usepackage{amsthm}
\usepackage{amssymb}
\usepackage{lastpage}
\usepackage{fancyhdr}
\usepackage{fancybox}
\usepackage{mathrsfs}

\def\cal{\mathcal}

\def\phi{\varphi }
\def\epsilon{\varepsilon}

\theoremstyle{plain}
\newtheorem{theorem}{Theorem}[section]
\newtheorem{corollary}[theorem]{Corollary}
\newtheorem{lemma}[theorem]{Lemma}
\newtheorem{proposition}[theorem]{Proposition}

\theoremstyle{definition}
\newtheorem{definition}[theorem]{Definition}

\newtheorem{example}[theorem]{Example}
\newtheorem{deffacts}[theorem]{Definitions and facts}
\numberwithin{equation}{section}

\begin{document}

\title{Positivity of Gibbs states on distance-regular graphs}
\author{Michael Voit}
\address{Fakult\"at Mathematik, Technische Universit\"at Dortmund,
          Vogelpothsweg 87,
          D-44221 Dortmund, Germany}
\email{michael.voit@math.tu-dortmund.de}


\date{\today}

\maketitle

\begin{abstract} We study criteria which ensure that Gibbs states (often also called generalized vacuum states) on distance-regular graphs
  are positive. Our main criterion assumes that  the graph can be embedded into a growing family of distance-regular graphs.
For the proof of the positivity we then use  polynomial hypergroup theory and
  translate this positivity into the problem whether for $x\in[-1,1]$ the function $n\mapsto x^n$ has a positive integral representation
  w.r.t.~the orthogonal polynomials associated with the graph.
We apply our criteria to several  examples. For Hamming graphs and the infinite distance-transitive graphs we obtain
  a complete description of the positive Gibbs states.
\end{abstract}

\section{Introduction}

It is well known that vacuum states on distance-regular graphs lead to interesting models in quantum probability;
see e.g. the monograph \cite{HO} of Hora and Obata and references there.
Besides these classical states one can also study Gibbs states on these graphs.
These states, which are also called generalized vacuum states in some
references like \cite{HO}, are related to Gibbs kernels on these graphs; see Section 2.3 of \cite{HO}. These Gibbs kernels on a 
distance-regular graph $\Gamma=(V,E)$ with vertex set $V$, edge set $E$, natural distance function $d$, and diameter
$D\in\mathbb N\cup\{\infty\}$ are defined by $Q_x(u,v):=x^{d(u,v)}$ ($u,v\in V$) for $x\in[-1,1]$ (with the convention $0^0=1$).
It can be easily seen that states associated with $Q_x$ are positive in the sense of quantum probability if and only if
$Q_x$ is positive semidefinite, i.e., if
 $$\sum_{i,j=1}^n c_i \bar{c_j} \> x^{d(u_i,u_j)}\ge0$$
for all $n\in\mathbb N$, $u_1,\ldots,u_n\in V$ and $c_1,\ldots,c_n\in \mathbb C$.
While this positivity is obvious for the vacuum kernel $Q_0$ as well for $Q_1$, the set
$$P_\Gamma:=\{x\in[-1,1]:\> Q_x \>\> \text{positive semidefinite}\}$$
is unknown for general distance-regular graphs.
On the other hand, in some simple cases like $D=1$, $P_\Gamma$ can be determined easily; see \cite{HO}.
Moreover, it was shown by Haagerup \cite{H} that $P_\Gamma=[-1,1]$ for the infinite homogeneous trees, and by the 
work of Bozejko \cite{Bo1, Bo2}, it was shown that $[0,1]\subset P_\Gamma$ holds for some classes of
distance-regular graphs like the Hamming and Johnson graphs; see also \cite{HO} for further details.

In this paper we present a further approach to prove the positivity
of some Gibbs kernels which gives some additional informations for further classes of examples.
The idea of the approach here is  as follow: For each  distance-regular graph $\Gamma$ there is a
canonical associated (usually finite) sequence of orthogonal polynomials $(P_k)_k$ with an orthogonality measure
with some support $S_\Gamma\subset[-1,1]$ where  $S_\Gamma$ is the spectrum of the graph.
A standard argument with a Bochner-type theorem for polynomial hypergroups (see \cite{BH, J, L} for the background)
now yields that a kernel $Q:\Gamma\times\Gamma\to\mathbb R$ of the form $Q(u,v)=f(d(u,v))$
for some function $f:\{0,1,...\}\to\mathbb R$
is  positive semidefinite if and only if there is a (necessarily unique)
positive measure $\mu$ on  $S_\Gamma$  with
$$Q(u,v)=f(d(u,v))=\check\mu(d(u,v))=\int_{S_\Gamma} P_{d(u,v)}(w)\> d\mu(w)\quad(u,v\in\Gamma).$$

We now assume that for a given distance-regular graph $\Gamma$ there is a sequence  $(\Gamma_n)_n$
of distance-regular graphs containing $\Gamma$ such that the coefficients of the three-term-recurrence relations of the associated
orthogonal polynomials converge (after a suitable normalization) to certain constants for $n\to\infty$,
which means that  the polynomials $P_k^n(x)$ associated with the graphs $\Gamma_n$ tend to $x^k$ for all finite $k\le D$ and $n\to\infty$.
We shall call this condition the infinite embedding property in Section 3.
With standard arguments on positive semidefiniteness it then follows that all accumulation points $x\in[-1,1]$ of the union
of the supports of all  $S_{\Gamma_n}$ are contained in $P_\Gamma$.

We shall see that this seemingly difficult criterion works quite well for several examples.
In particular, we obtain a precise description of $P_\Gamma$ in this way for all Hamming graphs in Section 5,
and we can also extend the result of  Haagerup \cite{H} to a precise description of $P_\Gamma$ for all known infinite, locally finite 
distance-regular graphs in Section 7. Furthermore, for the Johnson graphs we reprove the known fact  $[0,1]\subset P_\Gamma$
in Section 5, and for the Grassmann graphs over the finite fields $\mathbb F_q$ ($q$ a prime power), which are sometimes also called
$q$-Johnson graphs,
we obtain
$$ \{q^{-j}: \> j\in\mathbb N_0\}\cup\{0\}\subset P_\Gamma \quad\text{with}\quad \mathbb N_0:=\{0,1,2,\ldots\}.$$
We expect that our criterion can be also applied to  further classes of finite distance-regular graphs which are discussed e.g.~in \cite{APM, BI, BCN, DKT, KOT}.

We also point out that the approach of this paper for  distance-regular graphs via the associated
orthogonal polynomials can be extended to examples of higher rank, i.e., objects like buildings or suitable classes
association schemes for which the analogues of associated spherical functions are multivariate orthogonal polynomials;
see \cite{BI, V4, V5,Zi} and references there for a further reading.
In this case, however, one first has to agree about canonical extensions of the notions of Gibbs states.

The paper is organized as follows. Section 2 contains some known background material on   distance-regular graphs and
the associated orthogonal polynomials and polynomial hypergroups. In Section 3 we then study the  infinite embedding property
and its consequences. Sections 4-7 then are devoted to several classes of examples.

The author would like to thank H. Tanaka for some useful remarks.

\section{Distance-regular graphs and associated polynomial hypergroups}

In this introductory section we recapitulate some notations and facts on finite and
infinite distance-regular graphs and the associated 
orthogonal polynomials and polynomial hypergroups. We also briefly discuss Gibbs states in this context.
The main sources are \cite{HO, KOT} for  Gibbs states on  distance-regular graphs, \cite{BH, J} for basics on hypergroups,
and \cite{V4,V5} for the connections between distance-regular graphs and the  associated
polynomial hypergroups, where these connections are discussed there for association schemes in a more general context.

We begin with distance-regular graphs:

\subsection{Distance-regular graphs} Consider an undirected, connected graph $\Gamma=(V,E)$ with an at most countable set $V$ of vertices and a set $E$ of edges. Assume that
  the graph has no loops and is locally finite, i.e., each vertex has only finitely many neighbors.
  Let $d:V\times V\to \mathbb N_0$ be the usual distance and $D:=\sup_{x,y\in V} d(x,y)\in \mathbb N_0\cup\{\infty\}$ the diameter of
  $\Gamma$. Let $X:=\{k\in\mathbb N_0: k\le D\}$ be the set of all possible distances on  $\Gamma$.

  The graph $\Gamma=(V,E)$ is called  distance-regular if for all $i,j,k\in X$ and $x,y\in V$
  with $d(x,y)=k$, the number of all vertices $z\in V$ with $d(x,z)=i$ and $d(y,z)=j$ is independent of the choice of $x,y$. Hence, for
 $d(x,y)=k$,
  the numbers
  $$p_{i,j}^k:=|\{z\in V: \> d(x,z)=i, d(y,z)=j\}|$$
 depend only on $i,j,k\in X$.

 For $k\in X$ we now consider the adjacency matrices $A_k$ with entries
$$(A_k)_{x,y}:= \left\{ \begin{array}{cc}
\displaystyle 1 &
 {\rm if}\>\> d(x,y)=k
\\     0&  {\rm otherwise.}\\     \end{array} \right. 
 $$
 In particular $A_0$ is the identity matrix, and all $A_k$ are locally finite, i.e., all rows and columns have only finitely many
 non-zero entries.
 Moreover, for all $i,j\in X$
 $$A_i\cdot A_j= \sum_{k=|i-j|}^{i+j} p_{i,j}^k \> A_k$$
 where we agree in such sums that  $k\in X$ holds.
 In this way, the $\mathbb C$-linear span $\cal A(\Gamma):=span(A_i: \> i\in X)$ becomes a commutative and associative algebra consisting of
 symmetric matrices where this algebra is generated by $A_1$. In particular, for each $k\in X$,  $A_k$ is a polynomial of degree $k$ in $A_1$.
 For the details see \cite{HO}.

We now fix some vertex $o\in V$. The associated Gibbs state (or generalized vacuum state)
 with parameter $q\in[-1,1]$ is  defined as the linear functional $\phi_q:\cal A(\Gamma)\to\mathbb C$ with
 $$\phi_q(B):=\sum_{x\in V} q^{d(x,o)} B_{x,o} \quad\text{for}\quad B\in\cal A(\Gamma)$$
 where we agree that $0^0:=1$. With this agreement we have $\phi_0(B)=B_{o,o}=B_{x,x}$ for $x\in V$, i.e., $\phi_0$
 is the classical vacuum state.  Please notice that it not clear (except for $q=0,1$) that these Gibbs states
 are states on $\cal A(\Gamma)$ in the sense of quantum probability, i.e., that the (not necessarily locally finite) matrices  $Q_q$ with
 $$(Q_q)_{x,y}:= q^{d(x,y)} \quad\quad (x,y\in V)$$
 are positive semidefinite, i.e., that
 for all $n\in\mathbb N$, $x_1,\ldots,x_n\in V$ and $c_1,\ldots,c_n\in \mathbb C$,
 $$\sum_{i,j=1}^n c_i \bar{c_j} \> q^{d(x_i,x_j)}\ge0.$$
Mainly by the work of Bozejko \cite{Bo1, Bo2} it is known that for many families of distance-regular graphs, the Gibbs states 
$\phi_q$ are in fact states for $q\in[0,1]$ and sometimes also for $q$ in larger compact subintervals of $[-1,1]$. This is the case, for instance, for
the infinite homogeneous trees and the Johnson and Hamming graphs. For these and many further classes of examples
we refer to \cite{HO}.

We now study this positivity  in the  more general context of ``invariant''  kernels $Q:V\times V\to\mathbb C$ 
on  distance-regular graphs where $Q(x,y)$ depends on $d(x,y)$ only like for the Gibbs states. For this we translate the positivity of
$Q$  into the problem whether the associated function $f:X\to\mathbb C$ with $Q(x,y)=f(d(x,y))$
is positive definite on $X$ in some hypergroup sense.

To explain this we recapitulate some facts from \cite{V5} for commutative association schemes and the associated commutative hypergroups where
we restrict our attention here to  distance-regular graphs and the associated polynomial hypergroups.

For this,
we  start with the adjacency matrices $A_k$ ($k\in X$) of a distance-regular graph $\Gamma$. As above we fix $o\in V$ and define
  $$\omega_k:= |\{x\in V: \> d(x,o)=k\}|=p_{k,k}^0<\infty \quad(k\in X)$$
  as the numbers of vertices with distance $k$ from $o$. We now form the  renormalized stochastic matrices
  $$\tilde A_k:=\frac{1}{\omega_k} A_k \quad(k\in X).$$
  Then, for $i,j\in X$,
  \begin{equation}\label{tilde-product}
    \tilde A_i\cdot \tilde A_j= \sum_{k=|i-j|}^{i+j} \tilde p_{i,j}^k \>\tilde A_k \quad\text{ with}\quad
    \tilde p_{i,j}^k =\frac{\omega_k}{\omega_i\omega_j}p_{i,j}^k\ge0
    \end{equation}
  where  obviously  $\sum_{k=|i-j|}^{i+j} \tilde p_{i,j}^k =1$ holds.
  To get a stochastic interpretation of this identity, we recapitulate the following identity,
  which is well-known for association schemes (see \cite{BI} in the finite case or Lemma 3.5
  in \cite{V4} in general):
  \begin{equation}\label{asso-ident}
  \text{ For all}\quad i,j,k\in X:\quad\omega_k p_{i,j}^k =\omega_i p_{k,j}^i.
  \end{equation}
  Therefore, for all $i,j,k\in X$,
 \begin{equation}\label{tilde-product2}
    \tilde p_{i,j}^k  =\frac{1}{\omega_j}p_{k,j}^i\ge0.
    \end{equation}
 This means that for fixed $i,j\in X$, the  $\tilde p_{i,j}^k$ form the distribution for the distance from $o$,
 when we first make a random jump of size $i$ from $o$, and jump then again in an independent way with size $j$.

  With these notations we now define a convolution $*$ on the Banach space $M_b(X)$ of all bounded signed measures on $X$.
  In fact, for point measures $\delta_i,\delta_j$ ($i,j\in X$) we put
  \begin{equation}\label{convo-graph}
    \delta_i*\delta_j:= \sum_{k=|i-j|}^{i+j} \tilde p_{i,j}^k \>\delta_k,
       \end{equation}
  and extend $*$ to  $M_b(X)$ in a unique bilinear,
  weakly continuous way (the latter is necessary for $X=\mathbb N_0$ only).
  In this way, $(M_b(X),*)$ becomes a commutative Banach-$*$-algebra with
  the involution $.^*$ with $\mu^*(B):=\overline{\mu(B)}$ ($B\subset X, \> \mu\in M_b(X)$).
  More precisely, $(X,+)$ becomes a polynomial hypergroup, i.e., a special commutative discrete hypergroup.
  For the convenience of the reader we here briefly recapitulate the definition of these objects and collect a few facts:

\begin{definition}
  Let $X\ne\emptyset$ be an at most countable discrete set, and $*$ a
  weakly continuous, commutative, associative, bilinear product on the Banach space $M_b(X)$
 of all bounded signed  measures on $X$ with the following properties:
\begin{enumerate}
\item[\rm{(1)}] For all $i,j\in X$,  $\delta_i*\delta_j$  is a probability measure on $D$ with finite support.
\item[\rm{(2)}] There exists a neutral element $e\in E$ with $\delta_i*\delta_e=
\delta_e*\delta_i=\delta_i$ for  $i\in X$. 
\item[\rm{(3)}] There is an involution $x\mapsto\bar x$ on $X$
  such that for all $i,j\in X$, $e\in \text{supp}\> (\delta_i*\delta_j)$ holds
  if and only if $i=\bar j$. 
\item[\rm{(4)}] If for $\mu\in M_b(X)$, $\mu^-$ denotes the
  image of $\mu$ under the involution, then
  $(\delta_i*\delta_j)^-= \delta_{\bar j}*\delta_{\bar i}$ for all $i,j\in X$.
\end{enumerate}
Then $(X,*)$ is called a commutative discrete hypergroup.
$(X,*)$ is called symmetric if the involution is the identity. 
\end{definition}

We  collect some facts  on  commutative discrete hypergroups from \cite{BH,J}

\begin{deffacts}\label{facts-hypergroups}
  Let $(X,*)$ be a commutative discrete hypergroup.
\begin{enumerate}
\item[\rm{(1)}]  The identity $e$ and the involution $.^-$ above
are  unique.
\item[\rm{(2)}]  $(M_b(X),*)$ is a commutative
  Banach-$*$-algebra with the involution $\mu\mapsto\mu^*$ with
$\mu^*(A):=\overline{\mu(A^-)}$   for  $A\subset X$.
\item[\rm{(3)}] $(X,*)$ admits a Haar measure $\omega\in M^+(X)$, i.e., a nontrivial positive measure
  $\omega=\sum_{i\in X} \omega(i)\delta_i$ which satisfies $\omega=\omega*\delta_j=\delta_j*\omega$ for all
  $j\in X$. This Haar measure is unique up to a multiplicative constant. It can be defined by
  \begin{equation}\label{normalized-haar}
    \omega(i):=\frac{1}{\delta_i*\delta_{\bar i}(\{e\})} \quad \text{for}\quad i\in X.\end{equation}
  From now on we  use this Haar measure on $(X,*)$.
\item[\rm{(4)}] Let $C(X)$ be the space of all $\mathbb C$-valued functions on $X$,
  and $C_b(X)$ the subspace of all bounded functions on $X$. For $f\in C(X)$, $i,j\in X$ we write
  $$f(i*j):=(\delta_i*\delta_j)(f):=\sum_{k\in X} f(k)\cdot (\delta_i*\delta_j)(\{k\}).$$
\item[\rm{(5)}] The spaces of all (bounded) non-trivial multiplicative continuous functions on $(X,*)$ are
$$\chi(X,*):=\{\alpha\in C(X):\>\> \alpha\not\equiv 0,\>\>
\alpha(i* j)=\alpha(i)\cdot \alpha(j) \>\>\text{for all}\>\>
i,j\in X\}$$
and $\chi_b(X,*):=\chi(X,*)\cap C_b(X)$. Moreover,
$$\hat X:=(X,*)^\wedge:= \{\alpha\in \chi_b(X,*):\>\> \alpha(\bar i)=\overline{\alpha(i)} 
 \>\>\text{for all}\>\>
i\in X\}$$
is the dual space of $(X,*)$. Its elements are called characters. In particular the constant function 
${\bf 1}$ is a character.

If all spaces carry the topology of pointwise
convergence, then $\chi_b(X,*)$ and $\hat X$ are compact.
\item[\rm{(6)}] For $f\in L^1(X):= L^1(X,\omega)$ and $\mu\in M_b(X)$, their Fourier(-Stieltjes) transforms
  are 
$$\hat f(\alpha):=\int_X f(i)\overline{\alpha(i)}\> d\omega(i),\quad
\hat \mu(\alpha):=\int_X \overline{\alpha(i)}\> d\mu(i) 
\quad (\alpha\in\hat X).$$
 $\hat f$ and  $\hat\mu$ are continuous on $\hat X$ with
$\|\hat f\|_\infty\le\|f\|_1$, $\|\hat \mu\|_\infty\le\|\mu\|_{TV}$.
\item[\rm{(7)}] There exists a unique positive measure $\pi\in M^+(\hat X)$, called  the Plancherel measure on $\hat X$,
  such that the  Fourier transform
$.^\wedge: L^1(X)\cap L^2(X) \to C(\hat X)\cap L^2(\hat X,\pi)$
is an $L^2$-isometry.
The Fourier  transform $.^\wedge$ can be extended uniquely to an isometric isomorphism between 
$  L^2(X)$ and $ L^2(\hat X,\pi)$.
Moreover, by (\ref{normalized-haar}),  $\pi$ is a probability measure.

Notice that, different from  abelian groups, the support
$S:=supp\> \pi$ of $\pi$ may be a proper closed subset of $\hat X$.
In this case, we even have ${\bf 1}\not\in S$ for all known examples.

However, if $(X,*)$ is finite, then $S=\hat X$ holds with $|\hat X|=|X|$. 
\item[\rm{(8)}] For $f\in L^1(\hat X,\pi)$, $\mu\in M_b(\hat X)$, their
  inverse
Fourier transforms are given by
$$\check f(i):=\int_S f(\alpha) {\alpha(i)}\> d\pi(\alpha),\quad
\check \mu(i):=\int_{\hat D} {\alpha(i)}\> d\mu(\alpha) 
\quad (i\in X)$$
with $\check f\in C_0(X)$ (i.e. $\check f$ disappears at $\infty$) , $\check\mu\in C_b(X)$
and 
$\|\check f\|_\infty\le\|f\|_1$, $\|\check \mu\|_\infty\le\|\mu\|_{TV}$.
\item[\rm{(9)}] A function $f\in C_b(X)$ is called positive definite 
on  $(X,*)$ if for all $n\in\mathbb N$, $x_1,\ldots, x_n\in X$ and
  $c_1,\ldots,c_n\in\mathbb C$,
$\sum_{k,l=1}^n c_k \bar c_l \cdot f(x_k*\bar x_l)\ge0.$
Obviously, all characters $\alpha\in\hat X$ are positive definite.
\end{enumerate}
\end{deffacts}

The following theorem of Bochner (see \cite{J}) will be essential for our paper:

\begin{theorem}\label{Bochner}
  A function $f\in C_b(X)$ is positive
  definite if and only if $f=\check \mu$ for some  positive measure $\mu\in M_b^+(\hat X)$.
  This measure is unique, and $\mu$ is a probability measure if  and only if $\check\mu(e)=1$ holds.
\end{theorem}

In the context of homogeneous trees (i.e., infinite hypergroups $(X,*)$), we  also need the following variant; see
\cite{V2}:

\begin{proposition}\label{Bochner-supp}
  Let  $f\in C_b(X)$ be a positive definite function. Then  $f=\check \mu$ for some $\mu\in M_b^+(\hat X)$
with $supp\>\mu\subset S$ if and only if  $f$ is the pointwise limit of positive definite functions on $X$ with finite support.
\end{proposition}

We next turn our attention to polynomial hypergroups:

\begin{definition}\label{pol-hy}
  Let $D\in \mathbb N\cup\{\infty\}$ and $X:=\{i\in\mathbb N_0:\> i\le D\}$.  A commutative discrete hypergroup is called a polynomial
  hypergroup of diameter $D$, if there are numbers $a_i,b_i,c_i\ge0$ with $a_i+b_i+c_i=1$ ($i\in X$) with
  $a_0=1$, $b_0=c_0=0$, $a_i>0$ for $i\in X$ with $i<D$, $a_D=0$ for $D<\infty$, such that
  $$\delta_1*\delta_i=a_i\delta_{i+1}+b_i\delta_{i}+ c_i\delta_{i-1} \quad(i\in X).$$
\end{definition}

If we compare this definition with the convolution (\ref{convo-graph})
above in the context of distance-regular graphs, we see that each distance-regular graph $\Gamma$ leads to a polynomial
hypergroup structure $(X,*)$.
In particular, by (\ref{tilde-product2}),
the parameters  $a_i,b_i,c_i\ge0$ from Definition \ref{pol-hy} here are given by
\begin{equation}\label{abc-general}
  a_i=\tilde p_{1,i}^{i+1}=  \frac{1}{\omega_{1}} p_{1,i+1}^{i}  ,\quad  b_i=\tilde p_{1,i}^{i}=\frac{1}{\omega_{1}} p_{1,i}^{i},\quad
  c_i=\tilde p_{1,i}^{i-1}=\frac{1}{\omega_{1}} p_{1,i-1}^{i}.
       \end{equation}

We collect some classical facts on polynomial hypergroups; see \cite{BH,L}.
Let $(X,*)$ be a polynomial hypergroup with diameter $D$. Define  the polynomials $(P_i)_{i\in X}$  recursively by
\begin{equation}\label{recurrence-pol}
P_0= 1, \>\> P_1(x)=x,  \quad P_1\cdot P_i= a_iP_{i+1}+ b_iP_i+c_iP_{i-1} \quad (0< i<D).
\end{equation}
For  $x\in\mathbb C$, the functions $i\mapsto P_i(x)$ then form the multiplicative functions on  $(X,*)$,
where we need some additional restriction in the finite case which we discuss below. We first consider the infinite case.

\subsection{Multiplicative functions on infinite  polynomial hypergroups}
In this case, the 3-term-recurrence (\ref{abc-general}) and Favard's theorem (see e.g.~\cite{C}) show that 
$(P_i)_{i\ge0}$ is a sequence of orthogonal polynomials. Moreover, if we define the functions $\alpha_x(i):=P_i(x)$
for $x\in\mathbb C$, $i\in X=\mathbb N_0$, then $\alpha_1={\bf 1}$,
$\chi(X,*)=\{\alpha_x: \> x\in \mathbb C\}$ and  
$$\hat X=\{\alpha_x: \> x\in \mathbb R \>\>\text{with}\>\>(P_i(x))_{i\ge0} \>\>\text{bounded}\}.$$
If we identify $\chi(X,*)$ and $\hat X$ with $\mathbb C$ and the corresponding subset respectively, then the
topology of pointwise convergence agrees with the usual topology. In particular, $\hat X$ may be regarded as a compact subset
of $[-1,1]$. We  remark that then the Plancherel measure $\pi$ from Section \ref{facts-hypergroups}(7)
is  the orthogonality measure of $(P_i)_{i\ge0}$. It is a probability measure whose support is contained in $\hat X$.

\subsection{Multiplicative functions on finite  polynomial hypergroups}
Here, the 3-term-recurrence (\ref{abc-general}) also leads to a sequence $(P_i)_{i=0,\ldots,D}$ of orthogonal polynomials.
If we test whether the functions $\alpha_x$ as above are multiplicative, we land up with the condition
 $P_1(x)\cdot P_D(x)=b_D P_D(x) +c_DP_{D-1}(x)$ which is in fact solved for precisely
$D+1$ different points
$$-1\le x_D<x_{D-1}<\ldots<x_0=1.$$
In this case we  identify 
$\chi(X,*)=\hat X$ with $\{x_0,\ldots,x_D\}$.
The  Plancherel measure $\pi$ then is the orthogonality measure; its support is  $\{x_0,\ldots,x_D\}$.

\medskip

We  now turn to the central positivity result of this section.
We state it in the finite and infinite case separately, as  the infinite case is more involved. It follows from 
Theorems \ref{Bochner} and \ref{Bochner-supp} respectively, and it is shown in Section 6 of \cite{V5}
in the context of association schemes.

\begin{theorem}\label{main-pos1}
  Let $\Gamma=(V,E)$ be a finite distance-regular graph with diameter $D$. Then for a function $f:X\to\mathbb R$
  the following statements are equivalent:
\begin{enumerate}
\item[\rm{(1)}] The kernel $Q_f:V\times V\to\mathbb R$ with $Q_f(u,v):=f(d(u,v))$ is positive semidefinite;
\item[\rm{(2)}] The function $f$ is positive definite on the polynomial hypergroup $(X,*)$ associated with $\Gamma$;
\item[\rm{(3)}] There is a (unique) positive measure $\mu$ on $\hat X=\{x_0,\ldots,x_D\}$ with $f=\check\mu$.
\end{enumerate}
  \end{theorem}

\begin{theorem}\label{main-pos2} Let $\Gamma=(V,E)$ be an infinite distance-regular graph with associated polynomial hypergroup $(X,*)$.  Let
  $f:X\to\mathbb R$ be a function and  $Q_f$ the associated kernel on $V$ as before. Then:
  \begin{enumerate}
  \item[\rm{(1)}] If $Q_f$   is positive semidefinite, then  $f$ is a bounded positive definite function on  $(X,*)$, and there is
 a (unique) positive measure $\mu$ on $\hat X$ with $f=\check\mu$.
\item[\rm{(2)}] If $f$ is the pointwise limit of positive definite functions on  $(X,*)$ with finite supports, or if
  $f=\check\mu$ for some positive measure $\mu$ on $\hat X$ with $supp\> \mu\subset supp\> \pi$, then  $Q_f$
  is positive semidefinite.
  \end{enumerate}
\end{theorem}

Please notice that part (2) in Theorem \ref{main-pos2} is weaker than a complete converse statement of (1). We shall see below that
the complete converse statement of (1)  is not correct for some examples; see Section \ref{infinite-gr}.
We also remark that infinite polynomial hypergroups have unbounded positive definite functions.

\begin{example}\label{complete-graph}
  For an integer $N\ge2$ let $\Gamma$ be the complete graph with $N$ vertices, i.e., all vertices in $V:=\{1,\ldots,N\}$ are neighbored.
  Here the convolution $*$ of the associated polynomial hypergroup $(X=\{0,1\}, *)$ satisfies
$$\delta_0*\delta_0=\delta_0 ,\>\> \delta_0*\delta_1= \delta_1*\delta_0=\delta_1, 
  \>\> \delta_1*\delta_1= \frac{1}{N-1} \delta_0+\frac{N-2}{N-1}  \delta_1.$$
  Now let $x\in[-1,1]$. Then the Gibbs kernel $Q_x$ on  $V$ with $Q_x(u,v)=x^{d(u,v)}$
  is positive semidefinite if and only if the function $f_x(i):=x^i$ ($i=0,1$) is positive definite on $(X, *)$,
  and this is the case if and only if the matrix
  $$\begin{pmatrix} 1 &x \\
x & \frac{1}{N-1}+\frac{N-2}{N-1}x 
\end{pmatrix} $$
 is positive semidefinite which is the case by an elementary calculus precisely for $x\in [-1/(N-1),1]$.
\end{example}

\begin{example}\label{octahedron}
  The 6 vertices of an octahedron with its 12 edges form  a distance-regular graph $\Gamma=(V,E)$
  of diameter $D=2$; see Example 2.17 in \cite{HO}. The convolution $*$ of the associated polynomial
  hypergroup $(\{0,1,2\}, *)$ with identity $0$ satisfies
  $$\delta_2*\delta_2=\delta_0, \>\> \delta_1*\delta_2=\delta_1,
  \>\> \delta_1*\delta_1= \frac{1}{4}\delta_0+\frac{1}{4}\delta_2+\frac{1}{2}\delta_1.$$
  Therefore, by Theorem \ref{main-pos1}, the Gibbs state $Q_x(u,v):=x^{d(u,v)}$ is positive semidefinite on $V$
  if and only if the matrix
  $$\begin{pmatrix} 1 &x & x^2\\
    x & \frac{1}{4}+\frac{1}{4}x^2+\frac{1}{2}x &x\\
x^2&x&1\end{pmatrix} $$
is positive semidefinite. A computation of the principal minors
shows that this is the case precisely for $x\in[-2+\sqrt 3, 1]$; see also Example 2.17 in \cite{HO}.
\end{example}

\section{The infinite embedding property and positive Gibbs states}

We start with the following simple definition:

\begin{definition}
  Let $\Gamma_n:=(V_n,E_n)$ with $n=1,2$ be distance-regular graphs.
  $\Gamma_1$ is called a subgraph of $\Gamma_2$ (for short, $\Gamma_1\subset\Gamma_2$), if $V_1\subset V_2$, and if the distance function $d_2$ on $V_2$ restricted to $V_1$ is the distance function on $V_1$. 
\end{definition}

Notice that this subgraph property implies that vertices in $ V_1$ are
  neighbored in $V_1$ if and only if they are in  $V_2$, but that the the converse statement usually does not hold.

The following is obvious for distance-regular graphs $\Gamma_n:=(V_n,E_n)$:

\begin{lemma}\label{basic-subgraph} Let  $\Gamma_1\subset\Gamma_2$.
\begin{enumerate}
  \item[\rm{(1)}] The diameters $D_1,D_2$ and the associated spaces $X_1,X_2$ from subsection 2.1 satisfy $D_1\le D_2$ and $X_1\subset X_2$.
\item[\rm{(2)}]  If $Q:V_2\times V_2\to\mathbb R$ is a positive semidefinite kernel, then
  its restriction $Q|_{V_1}:V_1\times V_1\to\mathbb R$ is also positive  semidefinite.
\end{enumerate}\end{lemma}

This fact together with Theorems \ref{main-pos1} and \ref{main-pos2} implies:

\begin{lemma} Let $\Gamma_n:=(V_n,E_n)$ with $n=1,2$ be distance-regular graphs with $\Gamma_1\subset\Gamma_2$.
  Let $f\in C_b(X_2)$ be a function on the polynomial hypergroup $(X_2,*_2)$ associated with $\Gamma_2$
  such that the associated kernel $Q_f:V_2\times V_2\to\mathbb R$ with $Q_f(u,v):=f(d_2(u,v))$ is  positive semidefinite.
  Then the restriction $f|_{X_1}$ of $f$ is positive definite on the polynomial hypergroup $(X_1,*_1)$ associated with $\Gamma_1$,
  and the associated kernel $Q_{f|_{X_1}}=Q_f|_{V_1}$ on $V_1$ is positive semidefinite.
 \end{lemma}

In particular, by Theorem \ref{main-pos1}:

\begin{proposition}\label{sub-pos1} Let $\Gamma_n:=(V_n,E_n)$ with $n=1,2$ be finite distance-regular graphs with $\Gamma_1\subset\Gamma_2$.
  Let $f\in C_b(X_2)$ be a  positive definite function on the polynomial hypergroup $(X_2,*_2)$ associated with $\Gamma_2$.
  Then $f|_{X_1}$  is positive definite on the polynomial hypergroup $(X_1,*_1)$ associated with $\Gamma_1$,
  and the associated kernel $Q_{f|_{X_1}}=Q_f|_{V_1}$ on $V_1$ is positive semidefinite.
\end{proposition}

Moreover, by Theorem \ref{main-pos2}:

\begin{proposition}\label{sub-pos2} Let $\Gamma_n:=(V_n,E_n)$ with $n=1,2$ be infinite distance-regular graphs with $\Gamma_1\subset\Gamma_2$.
  Let $f\in C_b(X_2)$ be a   function on the polynomial hypergroup $(X_2,*_2)$ associated with $\Gamma_2$ such that
  $f$ is the pointwise limit of finitely supported positive definite functions on  $(X_2,*_2)$, or that
  $f=\check\mu$ for some positive measure $\mu$ on $\hat X_2$ with $supp\> \mu\subset supp\> \pi_2$ (with the Plancherel measure $\pi_2$
  on the dual  $\hat X_2$ of  $(X_2,*_2)$).
Then $f|_{X_1}$  is positive definite on the polynomial hypergroup $(X_1,*_1)$ associated with $\Gamma_1$,
  and the associated kernel $Q_{f|_{X_1}}=Q_f|_{V_1}$ on $V_1$ is positive semidefinite.
  \end{proposition}

The following definition is central for this section:

\begin{definition}\label{embedding} We say that a distance-regular graph  $\Gamma:=(V,E)$ with diameter $D$
  has the infinite embedding property if there is a sequence of distance-regular graphs $(\Gamma_n:=(V_n,E_n))_{n\in\mathbb N}$
  with $\Gamma\subset \Gamma_n$ for $n\in\mathbb N$ with the following property:

  Let $(X_n,*_n)$ be the polynomial hypergroups associated with the graphs $\Gamma_n$, and let
  $a_i^{(n)}, b_i^{(n)},c_i^{(n)}\ge0$ with $a_i^{(n)}+ b_i^{(n)}+c_i^{(n)}=1$  be the associated
              3-term recurrence coefficients as in (\ref{recurrence-pol}) for $i<D$. These coefficients satisfy
\begin{equation}\label{limit-recurrence}
\lim_{n\to\infty}a_i^{(n)}=1, \quad \lim_{n\to\infty}b_i^{(n)}=\lim_{n\to\infty}c_i^{(n)}=0 \quad(i<D).
       \end{equation} 
  \end{definition}

We notice, that by (\ref{abc-general}), the condition
(\ref{limit-recurrence}) means on the level of the  distance-regular graphs $\Gamma_n$
with the numbers $p_{i,j}^{k,(n)}, \omega_i^{(n)}$ from subsection
2.1 that
\begin{equation}\label{limit-recurrence-alt}
\lim_{n\to\infty} \frac{p_{1,i+1}^{i,(n)}}{\omega_1^{(n)}}=1 \quad(i<D).
\end{equation}

We  show in the next sections
that this infinite embedding property holds for some classical series of distance-regular graphs.
We here next discuss some consequences from this property.
We first consider the finite case:

\begin{theorem}\label{set-P}
  Let $\Gamma:=(V,E)$ be a finite distance-regular graph with the infinite embedding property with a corresponding
  sequence of distance-regular graphs $(\Gamma_n:=(V_n,E_n))_{n\in\mathbb N}$ and  associated  polynomial hypergroups $(X_n,*_n)$
  with the associated dual spaces  $\hat X_n\subset [-1,1]$. Let
  $$P:=\{x\in[-1,1]: \> x \quad\text{is an accumulation point of}\quad \bigcup_{n\in\mathbb N}\hat X_n\}.$$
  Then for all $x\in P$, the function $f_x(i):=x^i$ is positive definite on the hypergroup $(X,*)$ associated with $\Gamma$,
  and the kernel
  $Q_{f_x}(u,v):=x^{d(u,v)}$ is positive semidefinite on  $\Gamma$.
\end{theorem}

\begin{proof} Let $x\in P$. Then there exists a subsequence $(n_l)_{l\in \mathbb N}\subset\mathbb N$ and a 
  sequence $(x_l)_{l\in \mathbb N}\subset [-1,1]$ with $\lim_{l\to\infty}x_l=x$ such that for all $l$, $x_l\in \hat X_{n_l}$ holds.
  For each $l$ consider the orthogonal polynomials $(P_i^{(n_l)})_i$ associated with  the graph $\Gamma_{n_l}$.
  As the character $f_{x_l}(i):=P_i^{(n_l)}(x_l)$ ($i\in  X_{n_l}$) is positive definite on $( X_{n_l}, *_{n_l})$, 
  Proposition \ref{sub-pos1} implies that  that the function $f_{x_l}$ is  positive definite on $( X, *)$.
  On the other hand,  $\lim_{l\to\infty}x_l=x$, and the limit (\ref{limit-recurrence}) for the 3-term recurrence imply that
  the functions  $f_{x_l}$ tend to the function  $f_x$ from the theorem on $X$. This yields that $f_x$ is also 
   positive definite on $( X, *)$ as claimed. The last statement follows from Theorem \ref{main-pos1}.
\end{proof}
  
In the infinite case, the situation is slightly more complicated as  Theorem \ref{main-pos2} is weaker than
Theorem \ref{main-pos1}. The arguments of the preceding theorem thus only lead to the following weaker result.

\begin{theorem}\label{infinite-partial-positive}
  Let $\Gamma:=(V,E)$ be an infinite distance-regular graph with the infinite embedding property with the corresponding graphs
  $(\Gamma_n:=(V_n,E_n))_{n\in\mathbb N}$, $(X_n,*_n)$,
  $\hat X_n\subset [-1,1]$, and $P\subset[-1,1]$ as in the preceding theorem. 
  Then for all $x\in P$, the function $f_x(i):=x^i$ is positive definite on the hypergroup $(X,*)$ associated with $\Gamma$.\end{theorem}

On the other hand, in Section \ref{infinite-gr} we shall study a class of infinite  distance-regular graphs for which
we shall determine all $x\in [-1,1]$, for which  the Gibbs kernels $Q_{f_x}$ are  positive semidefinite.
To our knowledge, this class covers all known  infinite  distance-regular graphs.

In the end of this section we turn to another question.
It will turn out for most examples of distance-regular graphs $\Gamma$ in the next sections that
the Gibbs kernels $Q_x:=Q_{f_x}$ are positive semidefinite for all $x\in[-1,1]$.
In this context the following result may be of interest,
which shows that a weaker result already always ensures that the $Q_x$ are positive semidefinite for all $x\in[-1,1]$.
We  restrict our attention here to the finite case, as the infinite case will be treated completely in Section \ref{infinite-gr}
without this result.

\begin{proposition} 
  Let $\Gamma$ be a finite distance-regular graph. Assume that there is a sequence $(x_n)_{n\in\mathbb N}\subset ]0,1[$ with
      $\lim_{n\to\infty} x_n=1$ such that the Gibbs kernels $Q_{x_n}$ on $\Gamma$ as above are  positive semidefinite for all $n$. Then
      $Q_x$ is positive semidefinite for each $x\in[0,1]$.
\end{proposition}

\begin{proof}
  The set of positive semidefinite  kernels on $\Gamma$ is closed under pointwise limits,
  nonnegative linear combinations, and  pointwise products.
  This, the power series of the exponential function,  and the assumption yield that for each $n$ and each $t\ge0$ the kernel
  $$\tilde Q_{t,n}(u,v):= exp\Bigl( \frac{t}{1-x_n} (Q_{x_n}(u,v)-1)\Bigr) = e^{-\frac{t}{1-x_n}} \cdot exp\Bigl( \frac{t}{1-x_n} (Q_{x_n}(u,v)\Bigr)$$
with $u,v\in\Gamma$ 
      is  positive semidefinite. As
      $$\lim_{n\to\infty}\tilde Q_{t,n}(u,v)=exp(-t\cdot d(u,v))=Q_{e^{-t}}(u,v),$$
      we obtain that $Q_x$ is positive semidefinite for each $x\in]0,1]$. Finally, the case $x=0$ is obvious.
\end{proof}

In the next sections we discuss several examples.

\section{Hamming graphs and Krawtchouk polynomials}

Let $N\ge2,D\ge1$ be integers. Let $V:=\{1,2,\ldots,N\}^D$ be equipped with the metric
$$d(u,v):=|\{i=1,\ldots,D:\> u_i\ne v_i\}| \quad\quad(u,v\in V).$$
It is well known and can be easily seen that the graph $(V,E)$ with $E=\{\{u,v\}\in V^2:\> d(u,v)=1\}$
is distance-regular with diameter $D$, and with
\begin{equation}
  \omega_i=p_{i,i}^0= \binom{D}{i}\cdot (N-1)^i, \quad p_{1,i+1}^i= (D-i)(N-1) \quad(i=0,\ldots,D);
\end{equation}
see e.g. Section 5.1 of \cite{HO}. We denote this graph as Hamming graph $H(D,N)$.

For the Hamming graph $H(D,N)$, the 3-term recurrence coefficients in (\ref{abc-general}) satisfy
\begin{equation}\label{hamm1}
  a_i=\tilde p_{1,i}^{i+1}=\frac{p_{1,i+1}^{i}}{\omega_1}=\frac{D-i}{D}\quad(i=0,\ldots,D).
  \end{equation}  
Furthermore, for fixed  $N,D,\tilde D$ with $\tilde D>D$ we may regard $H(D,N)$ as subgraph of $H(\tilde D,N)$ in a canonical way.
Therefore, by (\ref{hamm1}) we see that  all  Hamming graphs $H(D,N)$ have the infinite embedding property
when we choose the graphs $\Gamma_n$ in Definition \ref{embedding} as $\Gamma_n:=H(D+n,N)$ for $n\in\mathbb N$.

We next identify the set $P$ in Theorem \ref{set-P}. For this we recall that the orthogonal polynomials $(P_i)_{i=0,\ldots,D}$ 
associated with the graphs  $H(D,N)$ are Krawtchouk polynomials up to affine-linear transformations.

For this we recapitulate e.g.~from Szeg\"{o} \cite{S} or Section 5 of \cite{DR} that for $p\in]0,1[$ 
 the Krawtchouk polynomials $K_i(x)=K_i(x;D,p)$ ($i=0,\ldots,D$) can 
be defined by
\begin{equation}
  K_i(x):=K_i(x;D,p):= \>_2F_1(-i,-x;-D;1/p)=\sum_{k=0}^N {(-i)_k(-x)_k\over (-D)_k k!}
  \Bigl(\frac{1}{p}\Bigr)^k
   \end{equation}  
By \cite{S, DR}, these polynomials have the following properties:
\begin{enumerate}
\item[\rm{(1)}] $K_l(x)=K_x(l)$ and   $K_0(x)=K_l(0)=1$ for   $x,l\in \{0,1,\ldots, D\}$;
\item[\rm{(2)}]   $K_1(x) = 1-x/(Dp)$;
\item[\rm{(3)}] Recurrence relation: For $l\in \{0,1,\ldots, D\}$,
$$K_1\cdot K_l = \frac{D -l}{D} K_{l+1} + \frac{2p-1}{ p} \frac{l}{ D}K_l +
\frac{1-p}{p} \frac{l}{ D}K_{l-1}.$$ 
\item[\rm{(4)}] Orthogonality: For $l,m\in  \{0,1,\ldots, D\}$,
$$\sum_{x=0}^D K_l(x) K_m(x) \cdot \binom{ D}{x} p^x (1-p)^{D-x} = 
 {D\choose l}^{-1} \Bigl({1-p\over p}\Bigr)^l \delta_{l,m}.$$
\end{enumerate}

We now put $p:=(N-1)/N$ for the  Krawtchouk polynomials and compare (2),(3) with
 the recurrence  (\ref{recurrence-pol})  for the
 $(P_i)_{i=0,\ldots,D}$ with the coefficients (\ref{hamm1}).
 This shows that
 $$K_i(x)= P_i\Bigl(1-\frac{Nx}{D(N-1)}\Bigr) \quad\quad i=0,\ldots,D.$$
In particular, the orthogonality measures of the polynomials  $(P_i)_{i=0,\ldots,D}$
have the supports
$$\Bigl\{1-\frac{Nx}{D(N-1)}:\> x=0,1,\ldots,D\Bigr\}.$$
We thus conclude that  the set $P$ in Theorem \ref{set-P} here is given by
$$P=[1-N/(N-1),1] =[-1/(N-1), 1].$$
Therefore,  Theorems \ref{set-P} and \ref{main-pos1} lead to $(1)\Longrightarrow (2)\Longleftrightarrow(3)$
in the following  result:

\begin{theorem}\label{set-P-Hamming}
  For a Hamming graph $H(D,N)$ and $x\in\mathbb R$ the following statements are equivalent:
\begin{enumerate}
\item[\rm{(1)}] $x\in [-1/(N-1), 1]$;
\item[\rm{(2)}]  The function $f_x(i):=x^i$ ($i\in X$) is positive definite on the hypergroup $(X,*)$ associated with
  $H(D,N)$;
\item[\rm{(3)}] The Gibbs kernel $Q_{f_x}(u,v):=x^{d(u,v)}$ is positive semidefinite on  $H(D,N)$.
\end{enumerate}
\end{theorem}

\begin{proof} Assume that (3) holds. As the graph $H(1,N)$ is just the complete
  graph of valency $N$ considered in Example
  \ref{complete-graph}, and as  $H(1,N)$ is a subgraph of $H(D,N)$, (1)
  follows from Lemma \ref{basic-subgraph} and the results in Example
  \ref{complete-graph}.

  All further conclusions were already shown above.
\end{proof}

\section{Johnson graphs}

We here mainly follow Section 6 of \cite{HO}.
Let $v\ge 1, D\ge0$ be integers with $D\le v/2$. Consider the Johnson graph $J(v,D)=(V,E)$ with
$$V:=\{x\subset\{1,\ldots,v\}: \> |x|=D\},\quad E:=\{\{x,y\}\subset V:\> |x\cap y|=d-1\}.$$

It is well known and can be checked by elementary combinatorics that  $J(v,D)$ is  distance-regular with diameter $D$
 with the parameters
\begin{equation}
  \omega_i=p_{i,i}^0= \binom{D}{i}\cdot \binom{v-D}{i}, \quad p_{1,i+1}^i= (D-i)(v-D-i) \quad(i=0,\ldots,D);
\end{equation}
see e.g. Lemmas 6.6 and 6.8 of \cite{HO}. 

For the Johnson graph $J(v,D)$ the  3-term recurrence coefficients $a_i$ in (\ref{abc-general}) thus satisfy
\begin{equation}\label{johnson1}
  a_i=\tilde p_{1,i}^{i+1}=\frac{p_{1,i+1}^{i}}{\omega_1}=\frac{(D-i)(v-D-i)}{D(v-D)}\quad(i=0,\ldots,D).
  \end{equation}  

We next show that the Johnson graphs have the infinite embedding property. For this we first observe:

\begin{lemma}\label{embedding-johnson}
Let $v,D$ as above. Then for integers $m,n\ge0$, the  Johnson graph $J(v,D)$ can be regarded as subgraph of $J(v+m+n, D+n)$.
\end{lemma}

\begin{proof} Clearly, $J(v,D)$ may be regarded as  subgraph of $J(v+m, D)$, when we only consider sets in the vertex set
  of $J(v+m, D)$ of size $D$, which only contain elements of $\{1,\ldots,v\}$.

  On the other hand, we also can regard  $J(v,D)$ as  subgraph of $J(v+n, D+n)$. For this we notice that 
  $J(v,D)$ is isomorphic with the Johnson graph $J(v,v-D)$ via taking the complement of a subset of size $D$ in 
  $\{1,\ldots,v\}$ where $J(v,v-D)$ is defined as above
  (in fact, the condition $D\le v/2$ was assumed only in order to avoid this doubling of the examples).
  In this way we obtain
  $$J(v,D)\sim J(v,v-D)\subset J(v+n,v-D)=J(v+n,(v+n)-(D+n))\sim J(v+n,D+n)$$
  which shows that we can regard  $J(v,D)$ as  subgraph of $J(v+n, D+n)$.

A combination of both arguments now completes the proof of the lemma.
\end{proof}

\begin{proposition}
  Let $v,D$ as above. Then the  Johnson graph $J(v,D)$ has the infinite embedding property with the sequence
  $(\Gamma_n:=J(v+2n,D+n))_n$ of graphs where  $J(v,D)\subset J(v+2n,D+n)$ holds for all $n\ge1$ as described
  in Lemma \ref{embedding-johnson}.
  
  Moreover, the associated set $P$ from Theorem \ref{set-P} is 
  $$ P=[0,1].$$
\end{proposition}

\begin{proof} By (\ref{johnson1}), the coefficients $a_i^{(n)}$ of the 3-term recurrence
  of the orthogonal polynomials associated with the graphs $J(v+2n,D+n)$ according to (\ref{abc-general}) satisfiy
  $$a_i^{(n)}=\tilde p_{1,i}^{i+1,(n)}=  \frac{p_{1,i+1}^{i,(n)} }{\omega_{1}^{(n)}} =
  \frac{(D+n-i)(v-D+n-i)}{(D+n)(v-D+n)} \longrightarrow 1$$
  for all $i=1,\ldots,D-1$ and $n\to\infty$. This proves the infinite embedding property.

  We next identify the set  $P$ from Theorem \ref{set-P}. For this we recapitulate from the literature
  (see Section 3.2 of \cite{BI} or Section 9.1 of \cite{BCN}) that the orthogonal polynomials associated
  the  Johnson graph $J(v,D)$ as described in Section 2 have a orthogonality measure which has the support
  $$S_{v,D}:=\Bigl\{ \frac{D(v-D)-j(v-j+1)}{D(v-D)}=1-\frac{j(v-j+1)}{D(v-D)}:\quad j=0,\ldots,D\Bigr\}.$$
  We now consider the numbers in this set for the Johnson graphs $J(v+2n,D+n)$. Then these numbers are given by
  $$x_{j,n}:= 1-\frac{j(v+2n-j+1)}{(D+n)(v-D+n)} \quad\quad (j=0,\ldots,D+n).$$
  As
  $$ x_{j,n}-x_{j+1,n}=\frac{v+2n-2j}{(D+n)(v-D+n)}\in [0, \frac{v+2n}{(D+n)(v-D+n)}],$$
  we see that the distances of the $ x_{j,n}$ tend to $0$  uniformly in $j$  for $n\to\infty$ with
  $x_{0,n}=1$ as largest and $x_{D+n,n}$ as smallest element where  $x_{D+n,n}\to 0$ holds. This yields $P=[0,1]$.
\end{proof}

As in the preceding section, this proposition and Theorem \ref{set-P} lead to the following result on Gibbs states
which is well known;
see e.g. Proposition 6.27 of \cite{HO} where this result is shown via the quadratic embedding test of Bozejko.

\begin{theorem}\label{positive-P-Johnson}
  Let $v,D$ as above and $x\in[0,1]$. Then the Gibbs kernel $Q_{f_x}(u,v):=x^{d(u,v)}$ is positive semidefinite on the
  Johnson graph $J(v,D)$.
\end{theorem}

Please notice that the interval $[0,1]$ is usually not optimal; see for instance Example \ref{octahedron} which is just the Johnson graph
 $J(2,4)$.

\section{$q$-analogues of Johnson graphs}

These graphs, which are also called Grassmann graphs in \cite{BCN}, are as follows: Let $S$ be a vector space of dimension $v\in\mathbb N$
over the finite field $\mathbb F_q$ with $q$ some power of a prime. 
Let $V$ be the finite set of all $D$-dimensional subspaces of $S$ with some positive integer $D\le v/2$.
Let
$$E:=\{\{x,y\}\subset V: \quad dim(x\cap y)=D-1\}.$$
It is well known (see e.g. Section III.6 of \cite{BI} or \cite{BCN}) that these graphs $J_q(v,D):=(V,E)$ are
distance-transitive with diameter $D$, and that
the distance of subspaces $x,y\in  V$ is given by
$$d(x,y)= D-dim(x\cap y).$$

The associated parameters in the sense of Section 2  can be also obtained from the literature. In particular, we have
\begin{equation}\label{johnson-q-1-start}
  \omega_1=p_{1,1}^0=\frac{(q^D-1)(q^{v-D}-1)}{(q-1)^2}\cdot q, \quad p_{1,i+1}^i= q^{2i+1}\cdot\frac{(q^{D-i}-1)(q^{v-D-i}-1)}{(q-1)^2};
\end{equation}
see e.g. Section III.6(ii) of \cite{BI}. 
For the $q$-Johnson graph $J_q(v,D)$ the  coefficients $a_i$ in the 3-term-recurrence (\ref{abc-general}) thus satisfy
\begin{equation}\label{johnson-q-1}
  a_i=\tilde p_{1,i}^{i+1}=\frac{p_{1,i+1}^{i}}{\omega_1}= \frac{(q^D-q^i)(q^{v-D}-q^i)}{ (q^D-1)(q^{v-D}-1)} \quad(i=0,\ldots,D).
  \end{equation}  

We next observe that the $q$-Johnson graphs also have the infinite embedding property like the usual Johnson graphs.
As the proof is completely analog to the preceding section, we skip the proof.

\begin{lemma}\label{embedding-q-johnson}
  Let $v,D,q$ as above. Then for integers $m,n\ge0$,
  the  $q$-Johnson graph $J_q(v,D)$ can be regarded as subgraph of $J_q(v+m+n, D+n)$.

  Moreover, $J_q(v,D)$ has the infinite embedding property with the sequence of graphs
  $(J_q(v+2n,D+n))_{n\ge1}$. 
\end{lemma}

Please notice that the last observation follows from the fact that the coefficients $a_i$ from (\ref{johnson-q-1})
for the graphs  $J_q(v+2n, D+n)$ tend to 1 for fixed $i$ for $n\to\infty$.

In order to determine the associated set  $P$ from Theorem \ref{set-P}, we recapitulate from Theorem 9.3.3 of \cite{BCN}
that the 
eigenvalues of the adjacency matrix $A_1$ for  $J_q(v,D)$  are given by
$$\theta_j= q^{j+1} \frac{(q^{D-j}-1)(q^{v-D-j}-1)}{(q-1)^2}- \frac{q^j-1}{q-1} \quad\quad(j=0,\ldots,D).$$
Taking our normalizations from Section 2 into account together with the equation (\ref{johnson-q-1-start})
for  $\omega_1$, we see that
the orthogonality measure of the orthogonal polynomials associated with $J_q(v,D)$ is supported by the points
\begin{align}\label{eigenvalues-q-Johnson}
  \frac{1}{(q^{D}-1)(q^{v-D}-1)}\Bigl(   (q^{D-j}-1) (q^{v-D-j}-1)q^j -(q^j-1)(q-1)/q\Bigr)\notag\\
  =
  \frac{1}{(q^{D}-1)(q^{v-D}-1)}\Bigl(q^{j-1} +q^{v-j}-q^{v-D}-q^D+(q-1)/q\Bigr)
  \end{align}
for $j=0,\ldots,D$. In particular, for $j=0$, these points are equal to 1, and for $j=D$ equal to
$$-\frac{q-1}{q(q^{v-D}-1)}.$$
We now consider the parameters $(v+2n, D+n)$ instead of $(v,D)$ for $n\to\infty$.
In this case the expression in the second line of
(\ref{eigenvalues-q-Johnson}) behaves like
$$ q^{-j}+ O(q^{-n}) \quad\quad(n\to\infty)$$
where the terms $O(q^{-n})$ are uniform in $j=0,\ldots,D$.
We thus obtain
$$P=\{q^{-j}: \> j\in\mathbb N_0\}\cup\{0\}.$$
This and Theorem \ref{set-P} lead to the following result for the  Gibbs kernels
which was stated in a more general setting and proved
by other methods recently  in Proposition 3.1 of \cite{KOT}.

\begin{theorem}\label{positive-P-Johnson-q}
  Let $v,D,q$ as above and $x\in \{q^{-j}: \> j\in\mathbb N_0\}\cup\{0\}$.
  Then the Gibbs kernel $Q_{x}(u,v):=x^{d(u,v)}$ is positive semidefinite on  $J_q(v,D)$.
\end{theorem}

This set $\{q^{-j}: \> j\in\mathbb N_0\}\cup\{0\}$ for the $q$-Johnson graphs is obviously not optimal.

\begin{example} Consider the $q$-Johnson graph $J_2(2,4)$. Here $D=2$, and the set $V$ of all 2-dimensional subspaces of $\mathbb F_2^4$
  has $15\cdot 14/6=35$ elements.
  Moreover, it can be easily checked that $\omega_0=1$, $\omega_1=18$, $\omega_2=16$.

  The convolution of point measures on $X=\{0,1,2\}$ with $0$ as identity has the form
  $$\delta_1*\delta_1=\frac{1}{\omega_1}\delta_0+\alpha_1\delta_1+ \beta_1\delta_2, \quad 
  \delta_2*\delta_2=\frac{1}{\omega_2}\delta_0+\beta_2\delta_1+ \alpha_2\delta_2, \quad
  \delta_1*\delta_2=\gamma_1\delta_1+\gamma_1\delta_2$$
  with suitable $\alpha_1,  \beta_1, \alpha_2,  \beta_2, \gamma_1, \gamma_2>0$.
  To compute these coefficients, we claim that the spaces
  $$A:=span((1,0,0,0), (0,1,0,0)), B:=span((1,0,0,0), (0,0,0,1))\in V$$
  with $d(A,B)=1$ admit 9 spaces $C\in V$ with  $d(A,C)=d(B,C)=1$.
  In fact there are 5 possibilities of spaces $C$ of the form
  $C=span((1,0,0,0), u)\in V$ with $u\in V\setminus(A\cup B)$ and 4 
  spaces $C$ of the form $C=span( u,v)\in V$ with $u\in A\setminus(A\cap B)$ and
  $v\in B\setminus(A\cap B)$. As $\delta_1*\delta_1$ is a probability measure, we conclude that
  \begin{equation}\label{q-j-example1}
    \delta_1*\delta_1=\frac{1}{18}(\delta_0+9\delta_1+ 8\delta_2).
    \end{equation}
  In principle, we can also determine the other convolution products in this combinatorial way.
However, the following approach may be more efficient.
   Our hypergroup $(X,*)$ fits into the class of hermitian hypergroups of order 3 considered in Section 4 of
   Wildberger \cite{W}. By p. 100 of \cite{W} we thus have
   $$\alpha_1= 1-\frac{1+\gamma_1\omega_2}{\omega_1}, \quad 
   \alpha_2= 1-\frac{1+\gamma_2\omega_1}{\omega_2}.$$
 These relations and (\ref{q-j-example1}) now yield $\gamma_1$, then $\gamma_2$, and finally  
 $ \alpha_2$ and $\beta_2$. In summary, we get
\begin{equation}\label{q-j-example2}
\delta_1*\delta_2=\frac{1}{2}(\delta_1+ \delta_2), \quad  \delta_2*\delta_2=\frac{1}{16}(\delta_0+9\delta_1+ 6\delta_2). 
\end{equation}
Therefore, by Theorem \ref{main-pos1}, the Gibbs kernel $Q_x(u,v):=x^{d(u,v)}$ is positive semidefinite on $V$
  if and only if the matrix
  $$D:=\begin{pmatrix} 1 &x & x^2\\
    x & \frac{1}{18}(1+9x+ 8x^2) &\frac{1}{2}(x+x^2)\\
x^2&\frac{1}{2}(x+x^2)&\frac{1}{16}(1+9x+ 6x^2) \end{pmatrix} $$
is positive semidefinite.
As
$$det \begin{pmatrix} 1 &x \\ x & \frac{1}{18}(1+9x+ 8x^2)\end{pmatrix}= \frac{1}{18}(1-x)(10x-1)$$
and $$det D=  \frac{1}{288}(1-x)^2(1-2x)(1+4x)(16x^2+18x+1),$$
and 
as the zeros of $16x^2+18x+1$ are approximately $x_1\sim -1,06$ and $x_2\sim -0,059$, we obtain
that  the Gibbs kernel $Q_x$ is positive if and only if
$x\in [x_2, 1/2]\cup\{1\}$ holds.
It is quite interesting that the general approach of the preceding sections and  Proposition 3.1 of \cite{KOT}``detect''
the upper bound $1/2$ of the interval.
This example is also interesting as it  is an example of a distance-regular graph where $Q_x$ is
not positive semidefinite for some 
$x\in [0,1]$.
\end{example}

\section{The infinite distance-transitive graphs}\label{infinite-gr}

The set of all infinite distance-transitive graphs 
is parametrized by two parameters as follows by
Macpherson \cite{Mp}.

 Let $a,b\ge2$ be integers, and $C_b$ the complete  graph
graph with $b$ vertices.
Consider the infinite
   graph $\Gamma:=\Gamma(a,b)$  where
 precisely $a$ copies of  $C_b$ are tacked together at each vertex 
in a tree-like way. For $b=2$, $\Gamma$ is the homogeneous tree of
valency $a$. 

We now consider the associated polynomial hypergroups $(X=\mathbb N_0,*)$.
Some counting shows (see \cite {V1a}) that the convolution $*$ satisfies
\begin{equation}\label{faltung}
\delta_m*\delta_n = \sum_{k=|m-n|}^{m+n} g_{m,n,k} \delta_k\in  M^1(\mathbb N_0)
\quad\quad  (m,n\in\mathbb N_0)
\end{equation}
with
$$ g_{m,n, m+n}= \frac{a-1}{a}>0, \quad
g_{m,n, |m-n|}=  \frac{1}{a(a-1)^{m\wedge n-1} (b-1)^{m\wedge n}}>0,$$
$$g_{m,n,|m-n|+2k+1}= \frac{b-2}{ a(a-1)^{m\wedge n-k-1}(b-1)^{m\wedge n-k}}\ge0
\quad  (k\le m\wedge n-1), $$
$$g_{m,n,|m-n|+2k+2}= \frac{a-2}{a(a-1)^{m\wedge n-k-1}(b-1)^{m\wedge n-k-1}}\ge0
\quad (k\le m\wedge n-2).$$
In particular,
\begin{equation}\label{3-term-rek-trees}
g_{n,1,n+1}=\frac{a-1}{a }, \quad 
g_{n,1,n}=  \frac{b-2}{a(b-1)}, \quad
g_{n,1,n-1}=\frac{1}{a(b-1)},
\end{equation}
We next define associated orthogonal polynomials $(P_n^{(a,b)})_{n\ge0}$
according to the general 3-term recurrence (\ref{recurrence-pol}) by
$$P_0^{(a,b)}:=1, \quad\quad 
 P_1^{(a,b)}(x):= x$$
and
\begin{equation}\label{recu}
P_1^{(a,b)}P_n^{(a,b)}= \frac{1}{a(b-1)}P_{n-1}^{(a,b)}
 + \frac{b-2}{a(b-1)} P_n^{(a,b)} +
\frac{a-1}{ a }P_{n+1}^{(a,b)} \quad\quad(n\ge1) .  
\end{equation}
These polynomials satisfy
\begin{equation}\label{prodcartier}\textstyle
 P_m^{(a,b)}P_n^{(a,b)}= \sum_{k=|m-n|}^{m+n} g_{m,n,k}P_k^{(a,b)} \quad (m,n\ge0). 
\end{equation}

Please notice that the polynomials $(P_n^{(a,b)})_{n\ge0}$ differ by some affine-linear transformation
from the corresponding  notations in \cite{V1a, V3, V4, V5, BH}.
More precisely, the polynomials
\begin{equation}\label{renorming-tree}
  \tilde P_n^{(a,b)}(x) :=  P_n^{(a,b)}  (T(x)) \quad\quad(n\in\mathbb N_0) \quad\text{with}
\quad T(x):=
  \frac{2}{a}\cdot\sqrt{\frac{a-1}{b-1}}\cdot x+
\frac{b-2}{a(b-1)}
\end{equation}
are the polynomials considered e.g. in \cite{V4, V5}.
Some formulas can be expressed more easily in terms of   $(\tilde P_n^{(a,b)})_{n\ge0}$.
For instance,
\begin{equation}\label{potenz}
 \tilde P_n^{(a,b)}\bigl(\frac{z+z^{-1}}{2}\bigr)= \frac{c(z)z^n +c(z^{-1})z^{-n}}
{((a-1)(b-1))^{n/2}} \quad\quad\text{for}\>\> z\in\mathbb C\setminus\{0, \pm 1\}
\end{equation}
with
\begin{equation}\textstyle
c(z):=\frac{(a-1)z -z^{-1} +(b-2)(a-1)^{1/2}(b-1)^{-1/2}}
{a(z-z^{-1})}.
\end{equation}
We  define
\begin{align}\label{xnull}
 \tilde s_0:=  \tilde s_0^{(a,b)}&:= \frac{2-a-b}{ 2\sqrt{(a-1)(b-1)}}, \quad  s_0:=T( \tilde s_0)=\frac{-1}{b-1}\notag\\
 \tilde s_1:= \tilde  s_1^{(a,b)}&:=\frac{ ab -a-b+2}{ 2\sqrt{(a-1)(b-1)}}, \quad s_1:=T( \tilde s_1)=1.
\end{align}
By \cite{V4}, the
 $ \tilde P_n^{(a,b)}$ fit into the Askey-Wilson
 scheme (pp.~26--28 of 
 \cite{AW}). In particular, by  \cite{AW}, the 
normalized orthogonality measure $ \tilde \rho= \tilde \rho^{(a,b)}\in M^1(\mathbb R)$
of the $ \tilde P_n^{(a,b)}$ is 
\begin{equation}\label{orthohne}
d \tilde \rho^{(a,b)}(x)= \tilde w^{(a,b)}(x)dx\Bigr|_{[-1,1]} 
\quad\quad{\rm for}\quad a\ge b\ge 2
\end{equation}
and 
\begin{equation}\label{orthmit}
d \tilde \rho^{(a,b)}(x)= \tilde w^{(a,b)}(x)
dx\Bigr|_{[-1,1]} + \frac{b-a}{b} d\delta_{\tilde s_0}
  \quad{\rm for}\quad b> a\ge 2
\end{equation}
with
$$ \tilde w^{(a,b)}(x):=\frac{a}{2\pi} \cdot  \frac{(1-x^2)^{1/2}}{( \tilde s_1-x)(x- \tilde s_0)}.$$

In summary, if we identify the dual space $\hat X$ of $(X,*)$ via the polynomials $ P_n^{(a,b)}$ as in Section 2
with a compact subset of $\mathbb R$, we
have the following observations:
\begin{enumerate}
\item[\rm{(1)}] $\hat X=[T(-\tilde s_1),1].$
\item[\rm{(2)}] The support $S$ of the Plancherel measure is equal to
  $$\Bigl[\frac{b-2}{a(b-1)}-\frac{2}{a}\cdot\sqrt{\frac{a-1}{b-1}}, \frac{b-2}{a(b-1)}+\frac{2}{a}\cdot\sqrt{\frac{a-1}{b-1}}\Bigr]$$
for $a\ge b\ge 2$, and  for $b>a\ge2$ the additional isolated point $s_0$ appears. 
\item[\rm{(3)}] $S=\hat X $ holds precisely for $a=b=2$.
 \item[\rm{(4)}] $s_0=T(-\tilde s_1)$ holds precisely if $a=2$ or $b=2$.
\end{enumerate}

The following theorem from \cite{V4}  is central for our considerations:

\begin{theorem}\label{podi} Let $x\in\mathbb R$. Then the kernel
$$\Gamma\times\Gamma\to\mathbb R, \quad
  (v_1,v_2)\longmapsto P_{d(v_1,v_2)}^{( a, b)}(x)$$
  is  positive semidefinite if and only if $x\in [s_0, 1]$ holds.
\end{theorem}

This
theorem has the following consequence which is interesting in view of the differences in the general
Theorems \ref{main-pos1} and \ref{main-pos2} in the finite and infinite case:

\begin{corollary}\label{podi-cor} Consider a graph $\Gamma:=\Gamma(a,b)$ with
  parameters $a,b\ge2$ as above, and let $(X,*)$ be the associated polynomial hypergroup.
  Then the following statements are equivalent:
  \begin{enumerate}
  \item[\rm{(1)}] $a=2$ or $b=2$;
\item[\rm{(2)}] If $f:X\to\mathbb R$ is a bounded positive definite function on  $(X,*)$, then 
  the associated kernel $Q_f:V\times V\to\mathbb R$ with $Q_f(u,v):=f(d(u,v))$ is positive semidefinite.
  \end{enumerate}
\end{corollary}
  
  \begin{proof} $(2)\Longrightarrow(1)$ is trivial by (4) above and Theorem \ref{podi}. On the other hand, as each
    bounded positive definite function on  $(X,*)$ has the form $f=\check\mu$ for some positive measure $\mu$
    on $\hat X $ by Bochner's theorem \ref{Bochner},  $(1)\Longrightarrow(2)$ also follows from (4) above and Theorem \ref{podi}.
  \end{proof}

  We are now ready for the main result of this section which generalizes a corresponding result of Haagerup \cite{H} for
  homogeneous trees.
  
\begin{theorem}\label{pos-kernel-tree}
  Consider a graph $\Gamma(a,b)$ with
  $a,b\ge2$ as above. Let $x\in\mathbb R$. Then the Gibbs kernel $Q_{x}(u,v):=x^{d(u,v)}$ is positive semidefinite on
$\Gamma(a,b)$ if and only if $x\in [\frac{-1}{b-1},1]$ holds.
\end{theorem}

\begin{proof}  It follows from  (\ref{3-term-rek-trees})
  that $\Gamma=\Gamma(a,b)$ has the infinite embedding property with
  the sequence $(\Gamma(a+n,b))_n$ of graphs with  $\Gamma\subset \Gamma(a+n,b)$. 
  As the parameter $ s_0=\frac{-1}{b-1}$ of the  graphs $\Gamma(a+n,b)$ does not depend on $n$, we see that the
 set $P$ in Theorem \ref{set-P} is
 $[\frac{-1}{b-1},1]$. Hence, by Theorem \ref{set-P},  $Q_{x}$ is positive semidefinite on
$\Gamma(a,b)$ for $x\in [\frac{-1}{b-1},1]$.

Now assume that $Q_{x}$ is positive semidefinite on
$\Gamma(a,b)$. Then $Q_{x}$ is also positive semidefinite on the complete  graph $C_b$ with $b$ vertices,
and the results in Example
  \ref{complete-graph} yield $x\in [\frac{-1}{b-1},1]$.
\end{proof}

In the end of the paper we briefly recapitulate some results on the spectral measures
$\mu_x\in M^1([-1,1])$ in the case of homogeneous trees (i.e. $b=2$) from Letac \cite{L} which seem to be unknown 
for a broader audience. We recapitulate that for fixed valency $a\ge2$ these measures are characterized via
\begin{equation}
x^n=\check\mu_x(n)= \int_{-1}^1 P_n^{(a,2)}(z)\> d\mu_x(z) \quad(n\in\mathbb N_0, x\in[-1,1]).
 \end{equation}
By p. 136 of \cite{L} we have
\begin{equation}\label{density-tree-general}
  d\mu_x(z)= \frac{a}{2\pi} \cdot  \frac{1-x^2}{1+(a-1)x^2-azx}\cdot
  \frac{\sqrt{\frac{4(a-1)}{a^2} -z^2}}{1-z^2}\Bigl|_{[-2\sqrt{a-1}/a, 2\sqrt{a-1}/a]} dz
 \end{equation}
for $|x|\le  \frac{1}{\sqrt{a-1}}$. Notice here that
$$1+(a-1)x^2-azx>0 \quad\quad\text{for}\quad |x|<  \frac{1}{\sqrt{a-1}}, \quad z\in [-2\sqrt{a-1}/a, 2\sqrt{a-1}/a],$$
and that  $1+(a-1)x^2-azx=0$ for $x=\pm \frac{1}{\sqrt{a-1}}$ and $z=({\rm sign}\> x)\cdot \frac{2\sqrt{a-1}}{a}$,
  in which case this denominator also appears as a factor in the square-root-part of (\ref{density-tree-general}), i.e., 
  the density of $\mu_x$ here has a singularity at one boundary point. Furthermore, it is mentioned on  p. 136 of \cite{L}
  that for $|x|\in ]1/\sqrt{a-1}, 1[$,  $\mu_x$ has also a density similar to (\ref{density-tree-general}) with one additional atom.
 For $x=\pm1$ we clearly have  $\mu_x=\delta_x$.

\end{document}